\documentclass[12pt]{article}

\usepackage{amsmath,amssymb,amsthm,amscd,a4wide,upref}

\usepackage[enableskew]{youngtab}
\usepackage{ytableau}

\usepackage{hyperref}
\hypersetup{colorlinks,citecolor=blue,filecolor=black,linkcolor=blue,urlcolor=blue}

\usepackage{lmodern}     
\usepackage[T1]{fontenc}

\begin{document}


\newcommand{\non}{\nonumber}
\newcommand{\scl}{\scriptstyle}
\newcommand{\sclnearrow}{{\scl\nearrow}\ts}
\newcommand{\scloplus}{{\scl\bigoplus}}
\newcommand{\wt}{\widetilde}
\newcommand{\wh}{\widehat}
\newcommand{\ot}{\otimes}
\newcommand{\fand}{\quad\text{and}\quad}
\newcommand{\Fand}{\qquad\text{and}\qquad}
\newcommand{\ts}{\,}
\newcommand{\tss}{\hspace{1pt}}
\newcommand{\lan}{\langle\ts}
\newcommand{\ran}{\ts\rangle}
\newcommand{\vl}{\tss|\tss}
\newcommand{\qin}{q^{-1}}
\newcommand{\tpr}{t^{\tss\prime}}
\newcommand{\spr}{s^{\tss\prime}}
\newcommand{\di}{\partial}
\newcommand{\hra}{\hookrightarrow}
\newcommand{\antiddots}
    {\underset{\displaystyle\cdot\quad\ }
    {\overset{\displaystyle\quad\ \cdot}{\cdot}}}
\newcommand{\dddots}
    {\underset{\displaystyle\quad\ \cdot}
    {\overset{\displaystyle\cdot\quad\ }{\cdot}}}
\newcommand{\atopn}[2]{\genfrac{}{}{0pt}{}{#1}{#2}}

\newcommand{\su}{s^{}}
\newcommand{\vac}{\mathbf{1}}
\newcommand{\vacf}{\tss|0\rangle}
\newcommand{\BL}{ {\overline L}}
\newcommand{\BD}{ {\overline D}}
\newcommand{\BE}{ {\overline E}}
\newcommand{\BP}{ {\overline P}}
\newcommand{\ol}{\overline}
\newcommand{\pr}{^{\tss\prime}}
\newcommand{\ba}{\bar{a}}
\newcommand{\bb}{\bar{b}}
\newcommand{\eb}{\bar{e}}
\newcommand{\bi}{\bar{\imath}}
\newcommand{\bj}{\bar{\jmath}}
\newcommand{\bk}{\bar{k}}
\newcommand{\bl}{\bar{l}}
\newcommand{\hb}{\mathbf{h}}
\newcommand{\gb}{\mathbf{g}}
\newcommand{\For}{\qquad\text{or}\qquad}
\newcommand{\OR}{\qquad\text{or}\qquad}
\newcommand{\emp}{\mbox{}}


\newcommand{\U}{{\rm U}}
\newcommand{\Z}{{\rm Z}}
\newcommand{\ZY}{{\rm ZY}}
\newcommand{\Ar}{{\rm A}}
\newcommand{\Br}{{\rm B}}
\newcommand{\Cr}{{\rm C}}
\newcommand{\Fr}{{\rm F}}
\newcommand{\Mr}{{\rm M}}
\newcommand{\Sr}{{\rm S}}
\newcommand{\Prm}{{\rm P}}
\newcommand{\Lr}{{\rm L}}
\newcommand{\Ir}{{\rm I}}
\newcommand{\Jr}{{\rm J}}
\newcommand{\Qr}{{\rm Q}}
\newcommand{\Rr}{{\rm R}}
\newcommand{\X}{{\rm X}}
\newcommand{\Y}{{\rm Y}}
\newcommand{\DY}{ {\rm DY}}
\newcommand{\Or}{{\rm O}}
\newcommand{\SO}{{\rm SO}}
\newcommand{\GL}{{\rm GL}}
\newcommand{\Spr}{{\rm Sp}}
\newcommand{\Zr}{{\rm Z}}
\newcommand{\ev}{{\rm ev}}
\newcommand{\Pf}{{\rm Pf}}
\newcommand{\Ann}{{\rm{Ann}\ts}}
\newcommand{\Norm}{{\rm Norm\tss}}
\newcommand{\Ad}{{\rm Ad}}
\newcommand{\SY}{{\rm SY}}
\newcommand{\Pff}{{\rm Pf}\tss}
\newcommand{\Hf}{{\rm Hf}\tss}
\newcommand{\trts}{{\rm tr}\ts}
\newcommand{\otr}{{\rm otr}}
\newcommand{\row}{{\rm row}}
\newcommand{\End}{{\rm{End}\ts}}
\newcommand{\Mat}{{\rm{Mat}}}
\newcommand{\Hom}{{\rm{Hom}}}
\newcommand{\id}{{\rm id}}
\newcommand{\middd}{{\rm mid}}
\newcommand{\ch}{{\rm{ch}\ts}}
\newcommand{\ind}{{\rm{ind}\ts}}
\newcommand{\Normts}{{\rm{Norm}\ts}}
\newcommand{\mult}{{\rm{mult}}}
\newcommand{\per}{{\rm per}\ts}
\newcommand{\sgn}{{\rm sgn}\ts}
\newcommand{\sign}{{\rm sign}\ts}
\newcommand{\qdet}{{\rm qdet}\ts}
\newcommand{\sdet}{{\rm sdet}\ts}
\newcommand{\Ber}{{\rm Ber}\ts}
\newcommand{\inv}{{\rm inv}\ts}
\newcommand{\inva}{{\rm inv}}
\newcommand{\grts}{{\rm gr}\ts}
\newcommand{\grpr}{{\rm gr}^{\tss\prime}\ts}
\newcommand{\degpr}{{\rm deg}^{\tss\prime}\tss}
\newcommand{\Cond}{ {\rm Cond}\tss}
\newcommand{\Fun}{{\rm{Fun}\ts}}
\newcommand{\Rep}{{\rm{Rep}\ts}}
\newcommand{\sh}{{\rm{sh}}}
\newcommand{\weight}{{\rm{wt}\ts}}
\newcommand{\chara}{{\rm{char}\ts}}
\newcommand{\diag}{ {\rm diag}}
\newcommand{\Bos}{ {\rm Bos}}
\newcommand{\Ferm}{ {\rm Ferm}}
\newcommand{\cdet}{ {\rm cdet}}
\newcommand{\rdet}{ {\rm rdet}}
\newcommand{\imm}{ {\rm imm}}
\newcommand{\ad}{ {\rm ad}}
\newcommand{\tr}{ {\rm tr}}
\newcommand{\gr}{ {\rm gr}\tss}
\newcommand{\str}{ {\rm str}}
\newcommand{\loc}{{\rm loc}}
\newcommand{\Gr}{{\rm G}}

\newcommand{\twobar}{{\bar 2}}
\newcommand{\threebar}{{\bar 3}}


\newcommand{\AAb}{\mathbb{A}\tss}
\newcommand{\CC}{\mathbb{C}}
\newcommand{\KK}{\mathbb{K}\tss}
\newcommand{\QQ}{\mathbb{Q}\tss}
\newcommand{\SSb}{\mathbb{S}\tss}
\newcommand{\TT}{\mathbb{T}\tss}
\newcommand{\ZZ}{\mathbb{Z}\tss}
\newcommand{\Sbb}{\mathbb{S}}
\newcommand{\ZZb}{\mathbb{Z}}


\newcommand{\Ac}{{\mathcal A}}
\newcommand{\Bc}{{\mathcal B}}
\newcommand{\Cc}{{\mathcal C}}
\newcommand{\Cl}{{\mathcal Cl}}
\newcommand{\Dc}{{\mathcal D}}
\newcommand{\Ec}{{\mathcal E}}
\newcommand{\Fc}{{\mathcal F}}
\newcommand{\Jc}{{\mathcal J}}
\newcommand{\Gc}{{\mathcal G}}
\newcommand{\Hc}{{\mathcal H}}
\newcommand{\Lc}{{\mathcal L}}
\newcommand{\Nc}{{\mathcal N}}
\newcommand{\Xc}{{\mathcal X}}
\newcommand{\Yc}{{\mathcal Y}}
\newcommand{\Oc}{{\mathcal O}}
\newcommand{\Pc}{{\mathcal P}}
\newcommand{\PD}{{\mathcal {PD}}}
\newcommand{\Qc}{{\mathcal Q}}
\newcommand{\Rc}{{\mathcal R}}
\newcommand{\Sc}{{\mathcal S}}
\newcommand{\Tc}{{\mathcal T}}
\newcommand{\Uc}{{\mathcal U}}
\newcommand{\Vc}{{\mathcal V}}
\newcommand{\Wc}{{\mathcal W}}
\newcommand{\Zc}{{\mathcal Z}}
\newcommand{\HC}{{\mathcal HC}}


\newcommand{\asf}{\mathsf a}
\newcommand{\bsf}{\mathsf b}
\newcommand{\csf}{\mathsf c}
\newcommand{\nsf}{\mathsf n}


\newcommand{\Sym}{\mathfrak S}
\newcommand{\h}{\mathfrak h}
\newcommand{\q}{\mathfrak q}
\newcommand{\n}{\mathfrak n}
\newcommand{\m}{\mathfrak m}
\newcommand{\p}{\mathfrak p}
\newcommand{\gl}{\mathfrak{gl}}
\newcommand{\oa}{\mathfrak{o}}
\newcommand{\spa}{\mathfrak{sp}}
\newcommand{\osp}{\mathfrak{osp}}
\newcommand{\g}{\mathfrak{g}}
\newcommand{\kgot}{\mathfrak{k}}
\newcommand{\agot}{\mathfrak{a}}
\newcommand{\bgot}{\mathfrak{b}}
\newcommand{\sll}{\mathfrak{sl}}
\newcommand{\f}{\mathfrak{f}}
\newcommand{\z}{\mathfrak{z}}
\newcommand{\Zgot}{\mathfrak{Z}}


\newcommand{\al}{\alpha}
\newcommand{\be}{\beta}
\newcommand{\ga}{\gamma}
\newcommand{\de}{\delta}
\newcommand{\De}{\Delta}
\newcommand{\Ga}{\Gamma}
\newcommand{\ep}{\epsilon}
\newcommand{\ee}{\epsilon^{}}
\newcommand{\ve}{\varepsilon}
\newcommand{\ls}{\ts\lambda\ts}
\newcommand{\vk}{\varkappa}
\newcommand{\vs}{\varsigma}
\newcommand{\vt}{\vartheta}
\newcommand{\ka}{\kappa}
\newcommand{\vp}{\varphi}
\newcommand{\la}{\lambda}
\newcommand{\La}{\Lambda}
\newcommand{\si}{\sigma}
\newcommand{\ze}{\zeta}
\newcommand{\om}{\omega}
\newcommand{\Om}{\Omega}
\newcommand{\up}{\upsilon}


\newtheorem{thm}{Theorem}[section]
\newtheorem{lemma}[thm]{Lemma}
\newtheorem{prop}[thm]{Proposition}
\newtheorem{cor}[thm]{Corollary}
\newtheorem{conj}[thm]{Conjecture}

\theoremstyle{definition}
\newtheorem{definition}[thm]{Definition}
\newtheorem{example}[thm]{Example}

\theoremstyle{remark}
\newtheorem{remark}[thm]{Remark}

\newcommand{\bth}{\begin{thm}}
\renewcommand{\eth}{\end{thm}}
\newcommand{\bpr}{\begin{prop}}
\newcommand{\epr}{\end{prop}}
\newcommand{\ble}{\begin{lemma}}
\newcommand{\ele}{\end{lemma}}
\newcommand{\bco}{\begin{cor}}
\newcommand{\eco}{\end{cor}}
\newcommand{\bex}{\begin{example}}
\newcommand{\eex}{\end{example}}
\newcommand{\bde}{\begin{definition}}
\newcommand{\ede}{\end{definition}}
\newcommand{\bre}{\begin{remark}}
\newcommand{\ere}{\end{remark}}
\newcommand{\bcj}{\begin{conj}}
\newcommand{\ecj}{\end{conj}}

\renewcommand{\theequation}{\arabic{section}.\arabic{equation}}

\numberwithin{equation}{section}


\newcommand{\bpf}{\begin{proof}}
\newcommand{\epf}{\end{proof}}


\def\beql#1{\begin{equation}\label{#1}}

\newcommand{\bal}{\begin{aligned}}
\newcommand{\eal}{\end{aligned}}
\newcommand{\beq}{\begin{equation}}
\newcommand{\eeq}{\end{equation}}
\newcommand{\ben}{\begin{equation*}}
\newcommand{\een}{\end{equation*}}

\title{\Large\bf Universal Capelli identities and quantum immanants for\\ the queer Lie superalgebra}

\author{{Iryna Kashuba\quad and\quad Alexander Molev}}

\date{} 
\maketitle


\begin{abstract}
We apply the recently introduced idempotents for the Sergeev superalgebra
to construct quantum immanants for the queer Lie superalgebra $\q^{}_N$
as central elements of its universal enveloping algebra.
We prove universal odd and even
Capelli identities for $\q^{}_N$ and use them to
calculate the images of the quantum immanants under the action of $\q^{}_N$ in differential operators.
We show that the Harish-Chandra images of the quantum immanants coincide
with the factorial Schur $Q$-polynomials.



\end{abstract}

\section{Introduction}
\label{sec:int}

By the classical Schur--Weyl duality,
the natural actions of the
symmetric group $\Sym_n$ and
the general linear Lie algebra $\gl_N$ on the space of tensors $(\CC^{N})^{\ot n}$
centralize each other. This implies a realization of polynomial $\gl_N$-modules via
the primitive idempotents $e_{\Uc}$ for $\Sym_n$ associated with the standard tableaux $\Uc$.
Namely, if $\Ec_{\Uc}$ is the operator in the tensor space corresponding to $e_{\Uc}$,
where the shape of $\Uc$ is a Young diagram $\la$ with at most $N$ rows, then the space
$\Ec_{\Uc}(\CC^{N})^{\ot n}$ affords an irreducible representation of $\gl_N$
whose highest weight is $\la$. As a consequence, by considering the character
of this representation, one recovers the trace expression for
the Schur polynomial $s_{\la}$,
\beql{trex}
\tr\ts \Ec_{\Uc} Y_1\dots Y_n=s_{\la}(y_1,\dots,y_N),
\eeq
where
\beql{tennot}
Y_a=1^{\ot(a-1)}\ot Y\ot 1^{\ot(n-a)}\in (\End\CC^N)^{\ot n}
\eeq
for the diagonal matrix $Y=\diag(y_1,\dots,y_N)$
and the trace is taken over all $n$ copies of the endomorphism algebra $\End\CC^N$.

The {\em Sergeev duality} originating in \cite{s:ta}
provides a similar pairing between the {\em Sergeev superalgebra} $\Sc_n$
and the queer Lie superalgebra $\q^{}_N$ acting on the superspace of tensors $(\CC^{N|N})^{\ot n}$.
Suitable idempotents $e_{\Uc}$ for $\Sc_n$ were introduced independently in the recent papers
\cite{kms:jm}
and \cite{ls:sb}.
They are parameterized by {\em standard barred tableaux} $\Uc$ whose shapes are shifted
Young diagrams $\la$. By taking the image $\Ec_{\Uc}$ of $e_{\Uc}$ acting in the space of tensors
we find that the space $\Ec_{\Uc}(\CC^{N|N})^{\ot n}$ is a representation of $\q^{}_N$
isomorphic to the direct sum of several copies of the simple modules associated with $\la$.
The counterpart of \eqref{trex} takes the form
\beql{strex}
\str\ts \Ec_{\Uc} Y_1\dots Y_n=Q_{\la}(y_1,\dots,y_N),
\eeq
where $Y$ now denotes
a diagonal matrix with the entries $y_1,\dots,y_N,-y_1,\dots,-y_N$ (see \eqref{y} below),
the supertrace is taken over all $n$ copies of the endomorphism algebra $\End\CC^{N|N}$
and $Q_{\la}(y)$ is the Schur $Q$-polynomial.
Formula \eqref{strex} is obtained as a consequence of the character evaluation
for the $\q^{}_N$-modules as we show below in Proposition~\ref{prop:char}.

Furthermore, we aim
to understand a noncommutative or `quantum' version of \eqref{strex}, where $Y$ is replaced by
the generator matrix $F$ for $\q^{}_N$. The corresponding element belongs to the center
$\Zr(\q^{}_N)$
of the universal enveloping algebra $\U(\q^{}_N)$ and we call it the {\em quantum immanant}
by analogy with Okounkov's work \cite{o:qi} for $\gl_N$. It takes the form
\beql{smuu}
\SSb_{\la}=\str\ts \Ec^{}_{\Uc}\ts\big(F_1+\ka_1(\Uc)\big)\dots \big(F_n+\ka_{\nsf}(\Uc)\big),
\eeq
where the {\em signed contents} $\ka_{\asf}(\Uc)$ are some constants depending on $\Uc$ as
defined in \eqref{sgncont} below.
We show that the image of $\SSb_{\la}$ under the Harish-Chandra isomorphism $\chi$ coincides with
a version of the
{\em factorial
Schur $Q$-polynomial} $Q^+_{\la}(y)$ introduced in \cite{i:ia}:
\beql{hchim}
\chi(\SSb_{\la})=Q^+_{\la}(y_1,\dots,y_N);
\eeq
see \eqref{qstar} below.
In particular, $\SSb_{\la}$ depends only on the shape of $\Uc$ but not on $\Uc$ itself.
Moreover, the quantum immanants $\SSb_{\la}$ with $\ell(\la)\leqslant N$
form a basis of
$\Zr(\q^{}_N)$ (Theorem~\ref{thm:hch}).

We thus generalize the original work of Nazarov~\cite{n:ci}, where
the {\em Capelli elements} $C_{\la}$ were introduced.
A principal result of \cite{n:ci}, which may be called
the {\em abstract Capelli identity}, is the equality of
the images of the centers of the universal enveloping algebras
of $\q^{}_N$ and $\q^{}_M$
in the context
of the $(\q^{}_N,\q^{}_M)$-Howe duality.
Moreover, explicit images of the Capelli elements
$C_{\la}$ in differential operators are described in Theorem~4.7 therein.
The elements $C_{\la}$ form a basis of the center $\Zr(\q^{}_N)$; their
connection with the quantum immanants $\SSb_{\la}$ is
discussed in Remark~\ref{rem:concz}.

A closely related basis $\{z_{\la}\}$ of $\Zr(\q^{}_N)$ was produced
by Alldridge, Sahi and Salmasian~\cite[Theorem~1.3]{ass:sq}
in relation with the {\em Capelli eigenvalue problem} for $\q^{}_N$.
The Harish-Chandra images $\chi(z_{\la})$ coincide with the polynomials $Q^+_{\la}(y)$,
up to scalar factors, so that the quantum immanant $\SSb_{\la}$
is proportional to $z_{\la}$; see \eqref{zlasla}. A connection between
the central elements $z_{\la}$ and $C_{\la}$ was established in
\cite[Corollary~1.5]{ass:sq} via their
Harish-Chandra images\footnote{The formulas for the Harish-Chandra images should be corrected
by changing the signs of the variables of the factorial
Schur $Q$-polynomials.\label{sahi}}; see \eqref{caim} below.

In our paper, we develop an alternative approach to devise underlying
{\em universal Capelli identities} for $\q^{}_N$.
Their proofs
turn out to be
simple direct calculations. As a consequence,
we get the images
of the quantum immanants $\SSb_{\la}$ in the differential operators (Corollary~\ref{cor:inde}).
We thus
obtain explicit formulas for the central elements $z_{\la}$
and the corresponding Capelli operators $D_{\la}$
of \cite{ass:sq}, and
generalize the Capelli identity of \cite{n:ci}.
This approach was inspired by a recent work of
Gurevich, Saponov and Zaitsev~\cite{gps:wt}, \cite{z:um},
where universal Capelli identities for the Lie algebra
$\gl_N$ and the quantized enveloping algebra $\U_q(\gl_N)$ were given.
The action of $\gl_N$
in the space of polynomials
in variables $x_{ia}$ with $i,a\in\{1,\dots,N\}$ is defined by
\ben
E_{ij}\mapsto \sum_{a=1}^N x_{ia}\di_{aj}
\een
for the standard basis elements $E_{ij}$, where we set $\di_{ai}=\di/\di x_{ia}$.
The identity which was pointed out in \cite{z:um} reads
\beql{ucagln}
(E_1-M^{(1)})\dots (E_n-M^{(n)})\mapsto X_1\dots X_nD_1\dots D_n,
\eeq
where we use notation \eqref{tennot} for the matrices $X=[x_{ia}]$, $D=[\di_{ai}]$
and $E=XD$, while $M^{(b)}$ is the image of the Jucys--Murphy element $(1,b)+\dots+(b-1,b)$
under the action of the symmetric group in the space of tensors $(\CC^N)^{\ot n}$.

Note that \eqref{ucagln} is, in fact, equivalent to the identities in \cite{o:yb}; their
simple proof
was later given in \cite{m:rh} and reproduced in \cite[Theorem~7.4.1]{m:yc}.
The latter identities can be expressed as
\beql{ucaidem}
(E_1-M^{(1)})\dots (E_n-M^{(n)})\ts\Ec_{\Uc}\mapsto X_1\dots X_nD_1\dots D_n\ts\Ec_{\Uc},
\eeq
for any standard tableau $\Uc$ with $n$ boxes.
Since the primitive idempotents
$e_{\Uc}$ form a decomposition of the identity, we recover \eqref{ucagln} by summing over all $\Uc$ in \eqref{ucaidem}\footnote{The same remark applies to the universal identities of \cite[Theorem~21]{gps:wt}
for $\U_q(\gl_N)$;
they are equivalent to the identities of \cite[Theorem~4.1]{jlm:ih}.}.
Furthermore,
by taking trace on both parts of \eqref{ucaidem} one gets the identities of \cite{o:qi}
which were also proved in \cite{n:yc},
as well as the earlier identities of \cite{w:so}; see also
recent work \cite{jlz:gc} for their generalizations and additional references.

Our universal odd Capelli identity (see Theorem~\ref{thm:uca})
is analogous to \eqref{ucagln}, but the role of the matrix algebras
is now played by the superalgebra $Q_N$. Namely, we have the identity
\beql{caqn}
(G_1+M^{(1)})\dots (G_n+M^{(n)})\mapsto X_1\dots X_nD_1\dots D_n,
\eeq
where $G\in Q_N\ot\U(\q^{}_N)$ is a generator matrix for $\q^{}_N$, while $X$ and $D$ are supersymmetric
counterparts of the matrices of variables and derivations, and
$M^{(b)}$ is the image of the {\em odd Jucys--Murphy element} $m_b\in\Sc_n$ (multiplied by $\sqrt{2}$)
under the action of the Sergeev superalgebra in the space of tensors $(\CC^{N|N})^{\ot n}$;
see Sec.~\ref{sec:uci} for definitions.

A similar even counterpart of \eqref{caqn} is proved in Theorem~\ref{thm:maf}.
It will be our main instrument to calculate the Harish-Chandra image
of the quantum immanant $\SSb_{\la}$ as given in \eqref{hchim}. Following
\cite[Prop.~4.8]{n:ci},
we use relation \eqref{strex} to show that the top degree component of the image
coincides with the Schur $Q$-polynomial $Q_{\la}(y)$. Then we use the annihilation
properties of the image and Ivanov's characterization theorem
for the factorial Schur $Q$-polynomials \cite{i:ia}
to conclude that the image coincides with $Q^+_{\la}(y)$.
This follows the same approach as one in \cite{o:qi}, \cite{o:yb} to calculate the Harish-Chandra
images of quantum immanants for $\gl_N$.

\section{Representations of the Sergeev superalgebra}
\label{sec:rsa}

The {\em Sergeev superalgebra} $\Sc_n$ \cite{s:ta} is the graded tensor product of two superalgebras
\ben
\Sc_n=\CC\Sym_n^-\ot \Cl_n,
\een
where $\CC\Sym_n^-$ is the {\em spin symmetric group algebra} generated by odd
elements $t_1,\dots,t_{n-1}$ subject to the relations
\ben
t_a^2=1,\qquad t_at_{a+1}t_a=t_{a+1}t_a t_{a+1},\qquad t_at_b=-t_bt_a,\qquad |a-b|>1,
\een
while
$\Cl_n$ is the {\em Clifford superalgebra} generated by odd elements $c_1,\dots,c_n$
subject to the relations
\ben
c_a^2=-1,\qquad c_ac_b=-c_bc_a,\qquad a\ne b.
\een
The group algebra $\CC\Sym_n$ of the symmetric group $\Sym_n$  is embedded in $\Sc_n$ so that
the adjacent transpositions $s_a=(a,a+1)\in\Sym_n$ are identified with the elements of $\Sc_n$ by
\ben
s_a=\frac{1}{\sqrt 2}\ts t_a(c_{a+1}-c_a).
\een
We thus get the alternative presentation of $\Sc_n$ as the semidirect product
$\CC\Sym_n\ltimes \Cl_n$ with the relations between elements of the symmetric group
and the Clifford superalgebra given by
\ben
s_ac_a=c_{a+1}s_a,\qquad s_ac_{a+1}=c_a s_a,\qquad s_ac_b=c_bs_a,\qquad b\ne a,a+1.
\een
The {\em Jucys--Murphy elements} $x_b\in\Sc_n$ were introduced in \cite{n:ys} and are given by
\beql{jmde}
x_1=0,\qquad x_b=\sum_{a=1}^{b-1}(a,b)(1+c_ac_b),\qquad b=2,\dots,n.
\eeq
Their odd counterparts $m_b$ were pointed out in \cite{s:hd}
with the use of the analogs of the transpositions in the Sergeev superalgebra
given by
\ben
t_{ab}=(-1)^{b-a-1}\ts t_{b-1}\dots t_{a+1}t_a t_{a+1}\dots t_{b-1},\qquad a<b,
\een
and $t_{ba}=-t_{ab}$. They have the form
\beql{oddjm}
m_1=0,\qquad m_b=t_{1b}+\dots+t_{b-1,b},\qquad b=2,\dots,n.
\eeq
The odd and even Jucys--Murphy elements are related by
\ben
x_a=\sqrt{2}\ts m_a\tss c_a,\qquad a=1,\dots,n,
\een
with $x_a^2=2\tss m_a^2$. Moreover, the $x_a$ pairwise commute, while the $m_a$
pairwise anticommute.

The irreducible representations of $\Sc_n$ are
parameterized by strict partitions
$\la=(\la_1,\dots,\la_{\ell})$ of $n$ with $\la_1>\dots>\la_{\ell}>0$
and $|\la|:=\la_1+\dots+\la_{\ell}=n$; see e.g. \cite[Ch.~3]{cw:dr}.
We will set $\ell(\la)=\ell$ to denote the length
of the partition and
write $\la\Vdash n$ to indicate that $\la$ is a strict partition on $n$.
A strict partition $\la$ is depicted by the shifted Young diagram obtained
from the usual Young diagram by shifting row $i$ to the right by $i-1$ unit boxes.
A ({\em shifted}) $\la$-{\em tableau} is obtained by writing the numbers $1,\dots,n$ bijectively
into the boxes of the shifted Young diagram $\la$.
Such a tableau
is called {\em standard} if its entries increase
from left to right in each row and from top to bottom in each column.
The entries in the boxes $(i,i)$ are called the {\em diagonal entries}.

The dimension
of the simple $\Sc_n$-module $U^{\la}$ is given by
\ben
\dim U^{\la}=2^{n-\lfloor \frac{\ell(\la)}{2}\rfloor}_{}\ts g_{\la},
\een
where $g_{\la}$ is the number of standard $\la$-tableaux
found by the Schur formula
\ben
g_{\la}=\frac{n!}{\la_1!\dots \la_{\ell}!}\ts \prod_{1\leqslant i<j\leqslant \ell}\ts
\frac{\la_i-\la_j}{\la_i+\la_j}.
\een

As in \cite{kms:jm}, we will be using {\em standard barred tableaux} obtained by
adding bars on non-diagonal entries:\footnote{Note that such tableaux
were introduced earlier in \cite{s:st} and \cite{w:ts}
in the context of shifted Robinson--Schensted correspondence. }
\ben
\begin{ytableau}
    1 & 2 & \ol{4} & \ol{5} & 8 & \ol{10}\\
    \none & 3 & \ol{6} & 9\\
    \none & \none & 7
\end{ytableau}
\ .
\een

We will use sans-serif letters $\asf, \bsf, \csf$, etc.
to denote arbitrary barred or unbarred
entries of $\Uc$.
Given a standard barred tableau $\Uc$, introduce
the {\em signed content} $\ka_{\asf}(\Uc)$
of any entry $\asf$ of $\Uc$
by the formula
\beql{sgncont}
\ka_{\asf}(\Uc)=\begin{cases}\phantom{-\ts}\sqrt{\si_a(\si_a+1)}
\qquad&\text{if $\asf=a$ is unbarred},\\[0.4em]
-\ts\sqrt{\si_a(\si_a+1)}\qquad&\text{if $\asf=\bar a$ is barred},
\end{cases}
\eeq
where $\si_a=j-i$ is the content of the box $(i,j)$ of $\la$
occupied by $a$ or $\bar a$. Similarly, we will
extend the notation for the Jucys--Murphy elements
by setting $x_{\asf}=x_{a}$ for $\asf=a$ and $\asf=\bar a$.

The following definition was given in \cite{kms:jm} (see also \cite{ls:sb} for an equivalent form).
For any standard barred tableau $\Uc$
introduce the element $e^{}_{\Uc}$ of the Sergeev superalgebra $\Sc_n$ by induction, setting
$e^{}_{\ts\tiny\young(1)}=1$ for the one-box tableau, and
\beql{murphyse}
e^{}_{\Uc}=e^{}_{\Vc}\ts \frac{(x_n-b_1)\dots
(x_n-b_p)}{(\ka-b_1)\dots (\ka-b_p)},
\eeq
where $\Vc$ is the barred tableau
obtained from $\Uc$ by removing the box $\al$ occupied by $n$ (resp. $\bar n$).
The shape of $\Vc$ is a shifted diagram which we denote by $\mu$ and
$b_1,\dots,b_p$ are the signed contents in all addable boxes of $\mu$ (barred and unbarred),
except for the entry $n$ (resp. $\bar n$), while $\ka$ is the signed content of
the entry $n$ (resp. $\bar n$).
According to \cite[Prop.~2.2]{kms:jm},
the elements $e^{}_{\Uc}$ are pairwise orthogonal idempotents in $\Sc_n$ which
form a decomposition of the identity.
Moreover,
\beql{xietsign}
x_{\asf}\ts e^{}_{\Uc}=e^{}_{\Uc}\ts x_{\asf}=\ka_{\asf}(\Uc)\ts e^{}_{\Uc},
\eeq
where $\ka_{\asf}(\Uc)$ is the signed content of the entry $\asf$ of\ \ts $\Uc$.

Recall the definition of the $\Sc_n$-modules $\wh U^{\la}$
associated with strict partitions $\la\Vdash n$ going
back to
\cite{n:ys}; see also \cite{hks:da} and \cite{ww:ls}.
Consider the vector space
\ben
\wh U^{\la}=\bigoplus_{\sh(\Tc)=\la}\ts \Cl_n \tss v^{}_{\Tc}
\een
with the basis vectors $v^{}_{\Tc}$ associated with the (unbarred) standard $\la$-tableaux $\Tc$.
The action of generators of $\Sc_n$ is given by
\begin{align}
s_a\tss v^{}_{\Tc}&=\Bigg(\frac{1}{\ka_{a+1}(\Tc)-\ka_{a}(\Tc)}
+\frac{c_ac_{a+1}}{\ka_{a+1}(\Tc)+\ka_{a}(\Tc)}\Bigg)\tss v^{}_{\Tc}
+\Yc_a(\Tc)\ts v^{}_{s_a\Tc},
\label{saavt}\\[0.5em]
x_a\tss v^{}_{\Tc}&=\ka_{a}(\Tc)\ts v^{}_{\Tc},
\non
\end{align}
for $a$ running over the sets $\{1,\dots,n-1\}$ and $\{1,\dots,n\}$, respectively,
where we assume that $v^{}_{s_a\Tc}=0$ if the tableau $s_a\Tc$ is not standard, and
\ben
\Yc_a(\Tc)=\sqrt{A\big(\ka_a(\Tc),\ka_{a+1}(\Tc)\big)},
\een
with
\ben
A(u,v)=1-\frac{1}{(u-v)^2}-\frac{1}{(u+v)^2}.
\een

The module $\wh U^{\la}$ is isomorphic to the direct sum of $2^{\lfloor \ell(\la)/2\rfloor}_{}$
copies of the simple module $U^{\la}$. The results of \cite{kms:jm} imply an explicit form
of such a decomposition which we will now describe; cf. \cite{ls:sb}.
Along with the basis vectors $v^{}_{\Tc}$ of $\wh U^{\la}$ consider the vectors
$v^{}_{\Uc}$ associated with the standard barred tableaux $\Uc$. Namely,
suppose that a certain tableau $\Tc$ is obtained from $\Uc$ by unbarring all
barred entries $\bar a_1,\dots,\bar a_r$ with $a_1<\dots< a_r$.
Then set
$
v^{}_{\ts\Uc}=c_{a_1}\dots c_{a_r}\ts v^{}_{\Tc}.
$
Observe that
$
x_{\asf}\ts v^{}_{\Uc}=\ka_{\asf}(\Uc)\ts v^{}_{\Uc}
$
for any entry $\asf$ of $\Uc$.
Hence the definition \eqref{murphyse} implies that
in $\wh U^{\la}$ we have
\beql{idema}
e^{}_{\Uc}\ts v^{}_{\Vc}=\de^{}_{\Uc\tss\Vc}\ts v^{}_{\Vc}
\eeq
for any two standard barred tableaux $\Uc$ and $\Vc$ with $n$ boxes.

If a standard barred tableau $\Uc$ of shape $\la$ has
diagonal entries $d_1<\dots<d_{\ell}$, then $e^{}_{\Uc}$ commutes
with the Clifford generators $c_{d_1},\dots,c_{d_{\ell}}$.
Consider the subalgebra $\Cl^{\ts\Uc}_{\ell}$ of $\Cl_n$ generated by $c_{d_1},\dots,c_{d_{\ell}}$.
Let $\de=(\de_1,\dots,\de_m)$ be an $m$-tuple with $\de_a\in\{1,-1\}$,
where $m$ is defined by $\ell=2m$ or $\ell=2m+1$ for the even and odd $\ell$,
respectively.
Introduce the corresponding
idempotent $\Ec^{\tss\Uc}_{\de}\in\Cl^{\ts\Uc}_{\ell}$ by
\beql{ideme}
\Ec^{\tss\Uc}_{\de}=\prod_{a=1}^m\frac{1+i\tss \de_a c_{d_{2a-1}}c_{d_{2a}}}{2}.
\eeq
Given $\de$,
denote by $U^{\la}_{\de}$ the $\Sc_n$-submodule of $\wh U^{\la}$, generated
by the vectors $\Ec^{\Uc}_{\de} v^{}_{\Uc}$, where $\Uc$ runs over all
standard barred tableaux $\Uc$ of shape $\la$.

\bpr\label{prop:decomp}
We have the direct sum decomposition of $\Sc_n$-modules
\ben
\wh U^{\la}=\bigoplus_{\de} U^{\la}_{\de},
\een
summed over all $m$-tuples $\de$. Moreover, each $\Sc_n$-module $U^{\la}_{\de}$ is irreducible
and isomorphic to $U^{\la}$.
\epr

\bpf
We start by verifying that the submodule $U^{\la}_{\de}$ is spanned over $\CC$ by the vectors
of the form
\beql{spancd}
c_{d_{j_1}}\dots c_{d_{j_s}} \Ec^{\tss\Uc}_{\de} v^{}_{\Uc},\qquad s\geqslant 0,
\eeq
where $1\leqslant j_1<\dots<j_s\leqslant \ell(\la)$, all $j_a$ are odd, and $\Uc$ runs over the
standard barred tableaux $\Uc$ of shape $\la$. It is clear that the span of
the vectors \eqref{spancd} is stable under the action of the Clifford superalgebra
$\Cl_n$. Hence, to check that the span is also stable under the action of
$\Sym_n\subset\Sc_n$, it is enough to verify that the application
of the generators $s_b\in\Sym_n$ to the vectors
with $s=0$ and all entries of $\Uc$ unbarred, produces vectors in the span.
It is immediate from \eqref{saavt} that if both $b$ and $b+1$ are
non-diagonal entries of $\Uc$, then the application of
$s_b$ to a vector in \eqref{spancd} equals a linear
combination of such vectors. Now suppose that $b=d^{}_{2r-1}$ for some $r$.
Then the entry $b+1$ of $\Uc$ has to be non-diagonal.
By \eqref{ideme} we have
\beql{sbcoe}
s_b\ts \Ec^{\tss\Uc}_{\de}=\prod_{a=1,\ a \ne r}^m
\frac{1+i\tss \de_a c_{d_{2a-1}}c_{d_{2a}}}{2}\ts \frac{1+i\tss \de_r c_{b+1}c_{e}}{2}\tss s_b,
\eeq
where we set $e=d^{}_{2r}$. Furthermore,
applying \eqref{saavt} again we find that
\beql{sbvt}
s_b\tss v^{}_{\Uc}=\frac{1}{\ka_{b+1}(\Uc)}\big(1+c_bc_{b+1}\big)\tss v^{}_{\Uc}
+\Yc_b(\Uc)\ts v^{}_{s_b\Uc}.
\eeq
Now write
\ben
(1+i\tss \de_r c_{b+1}c_{e})(1+c_bc_{b+1})=1+i\tss \de_r\ts c_{b}c_{e}
-i\tss \de_r\tss c_e\ts(1+i\tss \de_r c_{b}c_{e})\ts c_{b+1}
\een
and combine \eqref{sbcoe} and \eqref{sbvt} to conclude that
the application of $s_b$ in the case under consideration yields a linear
combination of vectors of the form \eqref{spancd}. The remaining cases
where $b+1=d^{}_{2r-1}$ and where $b$ or $b+1$ equals $d^{}_{2r}$ for some $r$,
are dealt with by a similar calculation.

Now observe that the inclusion
\beql{incl}
\wh U^{\la}\subseteq \sum_{\de} U^{\la}_{\de}
\eeq
holds due to the relation
\ben
v^{}_{\Uc}=\sum_{\de}\Ec^{\tss\Uc}_{\de} v^{}_{\Uc}
\een
for any standard barred tableau $\Uc$ of shape $\la$.
It remains to note that for any given $\de$, the number of vectors of the form \eqref{spancd}
coincides with the dimension of $U^{\la}$, while the sum of these numbers over all $\de$
equals $2^n\tss g_{\la}$ which coincides with $\dim\wh U^{\la}$.
This implies that the inclusion in \eqref{incl} is the equality and the sum is direct.
\epf

We point out the
following corollary implied by the above proof.

\bco\label{cor:bas}
Let $\de=(\de_1,\dots,\de_m)$ be an $m$-tuple with $\de_a\in\{1,-1\}$.
The vectors \eqref{spancd} form a $\CC$-basis of the simple submodule $U^{\la}_{\de}$
of $\wh U^{\la}$.
\qed
\eco

Note that explicit formulas for the
action of generators of the Sergeev superalgebra $\Sc_n$ in the basis of $U^{\la}_{\de}$ with
$\de=(1,\dots,1)$ were given in \cite[Corollary~3.4]{kms:jm}, where
an equivalent realization of
the module $U^{\la}$ as a left ideal of $\Sc_n$
was used.

It is well-known that the simple $\Sc_n$-module $U^{\la}$ is of type
$\mathsf{M}$ if $\ell(\la)$ is even. This means that $U^{\la}$ remains irreducible
over $\Sc_n$, regarded as an algebra (ignoring the $\ZZ_2$-grading).
If $\ell(\la)$ is odd, then $U^{\la}$ is
of type $\mathsf{Q}$ and so splits into the direct sum $U^{\la}=U^{\la+}\oplus U^{\la-}$
of two non-isomorphic irreducible modules over the algebra $\Sc_n$; see e.g \cite[Ch.~3]{cw:dr}.

\bco\label{cor:deul}
Suppose that $\ell=2m+1$ and let $\de=(\de_1,\dots,\de_m)$ with $\de_a\in\{1,-1\}$.
The module $U^{\la}_{\de}$ over the algebra $\Sc_n$ splits into the direct sum
\beql{ulapm}
U^{\la}_{\de}=U^{\la+}_{\de}\oplus U^{\la-}_{\de}
\eeq
of two irreducible submodules, where the basis vectors of $U^{\la\pm}_{\de}$ over $\CC$
are given by
\beql{spancdpm}
c_{d_{j_1}}\dots c_{d_{j_s}} \Ec^{\tss\Uc}_{\de}\ts \frac{1\pm(-1)^r i\ts c_{d_{\ell}}}{2}\ts
v^{}_{\Uc},\qquad s\geqslant 0,
\eeq
with $1\leqslant j_1<\dots<j_s\leqslant 2m-1$, all $j_a$ are odd, $\Uc$ runs over the
standard barred tableaux of shape $\la$ and $r=r(\Uc)$ is the number of barred
entries in $\Uc$.
\eco

\bpf
First we verify that each span of vectors \eqref{spancdpm} with a chosen $+$ or $-$,
is stable under the action
of $\Sc_n$. This is clear for the action of the Clifford algebra $\Cl_n$.
Arguing as in the proof of Proposition~\ref{prop:decomp},
we are left to consider the action of $s_b$ with $b=d_{\ell}$ and $b=d_{\ell}-1$
on the vectors
with $s=r=0$.
In the first case, we have
\ben
s_b\tss (1\pm i\ts c_{b}) = (1\pm i\ts c_{b+1})\ts s_b.
\een
Now use \eqref{sbvt} and write
\ben
(1\pm i\ts c_{b+1})(1+c_bc_{b+1})=1\pm i\ts c_{b}\pm i\tss (1\mp i\ts c_{b})\tss c_{b+1},
\een
confirming that the resulting vector belongs to
the span \eqref{spancdpm}.
A similar calculation applies in the case with $b=d_{\ell}-1$.

The proof is completed by noting that the number of vectors in each span in \eqref{spancdpm}
is half the dimension of the module $U^{\la}_{\de}$, while the latter
is contained in the sum of
$U^{\la+}_{\de}$ and $U^{\la-}_{\de}$.
\epf

It is well-known that representations of the Sergeev algebra $\Sc_n$ can also be regarded
as those of the finite group $\wt B_n$ (the {\em twisted hyperoctahedral group}),
where the central element $z\in \wt B_n$ acts as $-1$; see e.g. \cite[Sec.~3.3]{cw:dr}.
Therefore, some properties of representations of $\Sc_n$, such as the orthogonality
of matrix elements and characters are implied by the finite group theory.
Keeping in mind this connection, it is natural to choose a specific basis of $\Sc_n$ suitable
for the application of the theory.
In the summation formulas below, $h$ will run over the basis
elements of the form
\beql{basele}
\si\tss c_1^{\ve_1}\dots c_n^{\ve_n},\qquad \text{with}\quad\si\in\Sym_n\fand\ve_a\in\{0,1\};
\eeq
cf.~\cite{n:ci}. In particular, all such elements are invertible in $\Sc_n$.

\bpr\label{prop:sumu}
Let $\Vc$ be a fixed standard barred tableau of shape $\la$. We have the identity
\beql{idensu}
\frac{n!\ts 2^{\ell(\la)}}{g_{\la}}\sum_{\sh(\Uc)=\la} e^{}_{\Uc}=\sum_h h\ts e_{\Vc}\tss h^{-1}.
\eeq
\epr

\bpf
Since the algebra $\Sc_n$ is semisimple,
its representation
\beql{uhat}
U=\bigoplus_{\la\Vdash n}  U^{\la}
\eeq
is faithful. Therefore, it suffices to verify that the elements
on both sides of \eqref{idensu} act in $U$ as the same operator.

Suppose first that $\ell(\la)=2m$ is even and use the basis provided
by Corollary~\ref{cor:bas} for the representation $U^{\la}=U^{\la}_{\de}$
with $\de=(1,\dots,1)$. Denote the basis vector in \eqref{spancd} by $v^{}_{\Uc^j}$, where
$j$ denotes the tuple $(j_1,\dots,j_s)$. Recall some properties of the
idempotents $e^{}_{\Uc}$
pointed out in \cite[Sec.~2]{kms:jm}. The idempotent $e^{}_{\Uc}$
commutes with the Clifford generators $c_{d_1},\dots,c_{d_{\ell}}$
corresponding to the diagonal entries $d_1,\dots,d_{\ell}$ of $\Uc$. Moreover,
if $a$ is an unbarred non-diagonal entry of $\Uc$, then
\ben
e^{}_{\Uc}\ts c_a=c_a\ts e^{}_{\Uc'}\Fand e^{}_{\Uc'}\ts c_a=c_a\ts e^{}_{\Uc},
\een
where the tableau $\Uc'$ is obtained from $\Uc$ by replacing $a$ with $\bar a$.
These properties together with \eqref{idema} imply that
the application of the element on the left hand side of \eqref{idensu}
to a basis vector $v^{}_{\Uc^j}\in U^{\la}$ yields
$\frac{n!\ts 2^{\ell(\la)}}{g_{\la}}\ts v^{}_{\Uc^{j}}$.

Introduce
the matrices $[A_{\Wc^i\Uc^j}(h)]$ representing the basis elements $h\in\Sc_n$ acting in $U^{\la}$ by
\beql{matha}
h\ts v^{}_{\Uc^j}=\sum_{\Wc^i}A_{\Wc^i\Uc^j}(h)\ts v^{}_{\Wc^i},
\eeq
summed over the labelled standard barred tableaux $\Wc^i$, where $i$ runs over the tuples
of odd numbers as defined in \eqref{spancd}.
Then the application of the element on the right hand side of
\eqref{idensu} to the vector $v^{}_{\Uc^j}$ gives
\ben
\bal
\sum_h h\ts e_{\Vc}\tss \sum_{\Wc^i}A_{\Wc^i\Uc^j}(h^{-1})\ts v^{}_{\Wc^i}
{}&=\sum_h h\ts \sum_{i}A_{\Vc^i\Uc^j}(h^{-1})\ts v^{}_{\Vc^i}\\
{}&=\sum_h \sum_{i}A_{\Vc^i\Uc^j}(h^{-1})\ts \sum_{\Wc^k}A_{\Wc^k\Vc^i}(h)\ts v^{}_{\Wc^k}.
\eal
\een
The coefficient of the basis vector $v^{}_{\Wc^k}$ in this expansion equals
\ben
\sum_{i}\sum_h A_{\Vc^i\Uc^j}(h^{-1})\ts A_{\Wc^k\Vc^i}(h).
\een
However, by the orthogonality property of matrix elements
in irreducible representations of finite groups,
the sum
\ben
\sum_h A_{\Vc^i\Uc^j}(h^{-1})\ts A_{\Wc^k\Vc^i}(h)
\een
is zero, unless $\Wc^k=\Uc^j$. In this case, the sum equals $\dim\Sc_n/\dim U^{\la}$
and summing over the tuples $i=(i_1,\dots,i_r)$ of odd numbers $1\leqslant i_1<\dots<i_r\leqslant 2m-1$,
we get
\ben
2^m \ts\frac{\dim\Sc_n}{\dim U^{\la}}=\frac{2^{n+m}\tss n!}{2^{n-m}\ts g_{\la}}=
\frac{n!\ts 2^{\ell(\la)}}{g_{\la}}
\een
thus verifying the identity in the case under consideration.

We will use a similar argument in the case $\ell(\la)=2m+1$,
and consider the module \eqref{uhat}, where we take
$U^{\la}=U^{\la}_{\de}$
with $\de=(1,\dots,1)$. However,
we need
to refine the decomposition in \eqref{uhat} by splitting $U^{\la}=U^{\la+}\oplus U^{\la-}$ as given
in \eqref{ulapm} for each $\la\Vdash n$
so that the components are simple
modules over the {\em algebra} $\Sc_n$. Denote the basis vectors
in \eqref{spancdpm} by $v^{}_{\Uc^{j+}}\in U^{\la+}$ and $v^{}_{\Uc^{j-}}\in U^{\la-}$, where
$j$ denotes the tuple $(j_1,\dots,j_s)$.

As with the previous case,
the application of the element on the left hand side of \eqref{idensu}
to a basis vector $v^{}_{\Uc^{j\pm}}\in U^{\la}$ yields
$\frac{n!\ts 2^{\ell(\la)}}{g_{\la}}\ts v^{}_{\Uc^{j\pm}}$.
On the other hand,
by introducing
the matrices of the operators $h$ acting in the modules $U^{\la\pm}$
similar to \eqref{matha} and using the
orthogonality of matrix elements, we conclude that
the element on the right hand side of
\eqref{idensu} acting on $v^{}_{\Uc^{j\pm}}$ gives the same vector multiplied by
the constant
\ben
2^m \ts\frac{\dim\Sc_n}{\dim U^{\la\pm}}=\frac{2^{n+m}\tss n!}{2^{n-m-1}\ts g_{\la}}=
\frac{n!\ts 2^{\ell(\la)}}{g_{\la}},
\een
thus completing the proof.
\epf

Now introduce the {\em character} $\chi^{\la}$ of the module $U^{\la}$
as the element
\ben
\chi^{\la}=\sum_h \chi^{\la}(h)\ts h^{-1}\in\Sc_n,
\een
summed over the basis elements $h$ of the form \eqref{basele}, where
$\chi^{\la}(h)$ is the trace of $h$ acting in $U^{\la}$. We have the
following expression for $\chi^{\la}$ implied by Proposition~\ref{prop:sumu};
cf. \cite[Cor.~3.5]{n:ci}.

\bco\label{cor:char}
For any
fixed standard barred tableau $\Vc$ of shape $\la$ we have the identity
\beql{chidensu}
2^{\lfloor \frac{\ell(\la)}{2}\rfloor}_{}\chi^{\la}=\sum_h h\ts e_{\Vc}\tss h^{-1}.
\eeq
\eco

\bpf
Suppose first that $\ell(\la)=2m$. By using the basis $\{v^{}_{\Uc^j}\}$ of $U^{\la}$
introduced in the proof of Proposition~\ref{prop:sumu}, we can write
\ben
\chi^{\la}(h)=\sum_{\Uc^j}A_{\Uc^j\Uc^j}(h).
\een
On the other hand, we have the formula
\beql{thipo}
e^{}_{\Uc}=\frac{\dim U^{\la}}{\dim\Sc_n}\ts\sum_j\sum_h\ts A_{\Uc^j\Uc^j}(h)\tss h^{-1},
\eeq
which is verified by the same argument as in the proof of Proposition~\ref{prop:sumu}: apply
the elements on both sides to the basis vectors of the module \eqref{uhat}
and use the orthogonality of matrix elements. Hence,
\ben
\sum_{\sh(\Uc)=\la} e^{}_{\Uc}=\frac{\dim U^{\la}}{\dim\Sc_n}\ts \chi^{\la},
\een
and \eqref{chidensu} follows from \eqref{idensu}.

Similarly, if $\ell(\la)=2m+1$, then using the bases $\{v^{}_{\Uc^{j\pm}}\}$ of $U^{\la\pm}$,
we have
\ben
\chi^{\la}(h)=\sum_{\Uc^{j+}}A_{\Uc^{j+}\Uc^{j+}}(h)+\sum_{\Uc^{j-}}A_{\Uc^{j-}\Uc^{j-}}(h),
\een
where we introduced the diagonal matrix elements $A_{\Uc^{j\pm}\Uc^{j\pm}}(h)$ of $h$
acting on the basis vectors. The counterpart of \eqref{thipo} takes the form
\beql{thipoodd}
e^{}_{\Uc}=\frac{\dim U^{\la+}}{\dim\Sc_n}\ts\sum_{j+}\sum_h\ts A_{\Uc^{j+}\Uc^{j+}}(h)\tss h^{-1}
+\frac{\dim U^{\la-}}{\dim\Sc_n}\ts\sum_{j-}\sum_h\ts A_{\Uc^{j-}\Uc^{j-}}(h)\tss h^{-1},
\eeq
which implies
\ben
\sum_{\sh(\Uc)=\la} e^{}_{\Uc}=\frac{\dim U^{\la}}{2\dim\Sc_n}\ts \chi^{\la}.
\een
Now \eqref{chidensu} follows from \eqref{idensu}.\footnote{A similar proof applies to
the well-known counterpart of \eqref{chidensu} for the characters of the symmetric
group.
Note that
the argument verifying it in \cite[Cor.~7.4.2]{m:yc} is erroneous.}
\epf

\section{The Schur $Q$-polynomials and Sergeev duality}
\label{sec:sqs}

We start by recalling the combinatorial
definition of the Schur $Q$-polynomials and their factorial counterparts.
Consider the ordered alphabet $1'<1<2'<2<\dots<N'<N$ and set
\ben
|k|=|k'|=k,\qquad \sgn k=-\sgn k'=1,
\een
for its entries. A {\em marked shifted Young tableau $\Tc$ of shape $\la$}
is obtained by labelling the boxes of $\la$ by the symbols of the alphabet in such a way that
the labels increase weakly along each row and down
each column. Moreover, for every $k\in\{1,\dots,N\}$, there is at most
one $k$ in each column and at most one $k'$ in
each row. The {\em factorial Schur $Q$-polynomials} $Q^+_{\la}(y)$ and $Q^-_{\la}(y)$
in the set of variables $y=(y_1,\dots,y_N)$ are given by
\beql{qstar}
Q^{\pm}_{\la}(y)=\sum_{\sh(\Tc)=\la}
\prod_{\al\in\la}\big(y^{}_{|\Tc(\al)|}\pm\sgn(\Tc(\al))\si(\al)\big),
\eeq
summed over the marked shifted Young tableau $\Tc$ of shape $\la$, the product
is taken over all boxes $\al$ of $\la$, and $\si(\al)=j-i$ denotes the content of $\al=(i,j)$.
Both $Q^+_{\la}(y)$ and $Q^-_{\la}(y)$
are particular cases of more general polynomials associated with arbitrary parameter sequences
$a$, as defined in \cite{i:ia}, where the polynomials $Q^-_{\la}(y)$ are denoted by $Q^*_{\la}(y)$.
Obviously, $Q^+_{\la}(y)=(-1)^{|\la|}\ts Q^-_{\la}(-y)$.
The top degree component of $Q^{\pm}_{\la}(y)$ with $\ell(\la)\leqslant N$ is the
{\em Schur $Q$-polynomial}
\ben
Q_{\la}(y)=\sum_{\sh(\Tc)=\la}\prod_{\al\in\la} y^{}_{|\Tc(\al)|};
\een
see e.g. \cite[Sec.~III.8]{m:sf}.

The algebra of {\em supersymmetric polynomials} $\Ga_N$ in the set of variables $y=(y_1,\dots,y_N)$
consists of those symmetric polynomials $P(y)$ in $y$ which satisfy the {\em cancellation property}:
the result of setting $y_1=-y_2=z$ in $P(y)$ does not depend on $z$.
The polynomials $Q_{\la}(y)$ and $Q^{\pm}_{\la}(y)$ are supersymmetric. Moreover, each of these
three families of polynomials corresponding to the shifted Young diagrams $\la$ with at most $N$ rows
forms a basis of $\Ga_N$.

By convention, any shifted Young diagram $\la$ with
$\ell(\la)\leqslant N$ will be regarded as the $N$-tuple 
$(\la_1,\dots,\la_N)=(\la_1,\dots,\la_{\ell},0,\dots,0)$.
The following {\em characterization property} holds:
if the top degree component of a polynomial $P(y)\in\Ga_N$ coincides with
$Q_{\la}(y)$ for some $\la$ with $\ell(\la)\leqslant N$, and $P(\mp\mu)=0$ for all
$\mu$ with $|\mu|<|\la|$, then $P(y)=Q^{\pm}_{\la}(y)$ \cite{i:ia}.

To recall the Sergeev duality,
consider the $\ZZ_2$-graded vector space $\CC^{N|N}$ over the field of complex numbers with the
canonical basis
$e_{-N},\dots,e_{-1},e_1,\dots,e_N$, where
the vector $e_i$ has the parity
$\bi\mod 2$ and
\ben
\bi=\begin{cases} 1\qquad\text{for}\quad i<0,\\
0\qquad\text{for}\quad i>0.
\end{cases}
\een
The endomorphism algebra $\End\CC^{N|N}$ is then equipped with a $\ZZ_2$-gradation with
the parity of the matrix unit $E_{ij}$ found by
$\bi+\bj\mod 2$. As in \cite{n:ci}, we also regard the $E_{ij}$ as basis elements
of the Lie superalgebra $\gl_{N|N}$. Furthermore, we consider the {\em queer Lie superalgebra}
$\q^{}_N$ as the subalgebra of $\gl_{N|N}$ spanned by the elements $F_{ij}=E_{ij}+E_{-i,-j}$.
This leads to a natural action of $\q^{}_N$ on the space $\CC^{N|N}$ and its tensor power
$(\CC^{N|N})^{\ot n}$. The action on the tensor product space commutes with the action
of the Sergeev superalgebra $\Sc_n$ defined by
\beql{saact}
c_a\mapsto J_a \Fand (a,b)\mapsto P_{ab},
\eeq
where we use a standard subscript notation for the copies of the element
\beql{j}
J=\sum_{i=-N}^N E_{i,-i}(-1)^{\bi}\in \End\CC^{N|N}
\eeq
and the
permutation operator
\ben
P=\sum_{i,j=-N}^N E_{ij}\ot E_{ji}(-1)^{\bj}\in \End\CC^{N|N}\ot\End\CC^{N|N},
\een
in the tensor product superalgebra $(\End\CC^{N|N})^{\ot n}$. Namely,
\beql{ja}
J_a=1^{\ot (a-1)}\ot J\ot 1^{\ot (n-a)},\qquad a=1,\dots,n,
\eeq
and
\ben
P_{ab}=\sum_{i,j=-N}^N 1^{\ot (a-1)}\ot E_{ij}\ot 1^{\ot (b-a-1)}\ot E_{ji}\ot 1^{\ot (n-b)}\ts(-1)^{\bj},
\een
for $1\leqslant a<b\leqslant n$. The zero values of the indices $i,j$, etc., are understood to be skipped
in such summation formulas. We will denote by $\Ec^{}_{\Uc}\in (\End\CC^{N|N})^{\ot n}$ the image
of the idempotent $e^{}_{\Uc}\in \Sc_n$ defined in \eqref{murphyse} under this action.

By the Sergeev duality originating in \cite{s:ta},
we have a decomposition of the $\U(\q^{}_N)\ot\Sc_n$-module
\ben
(\CC^{N|N})^{\ot n}\cong \bigoplus_{\la\Vdash n,\ \ell(\la)\leqslant N}\ts 2^{-\de(\la)}L(\la)\ot U^{\la},
\een
where $L(\la)$ is the simple $\q^{}_N$-module with the highest weight $\la$, while
$\de(\la)=0$ or $1$ depending on whether $\ell(\la)$ is even or odd. In the case of odd $\ell(\la)$,
the tensor product $L(\la)\ot U^{\la}$ is a direct sum of two isomorphic copies
of a simple module, and the factor $2^{-1}$ indicates that only one of
the copies occurs in the decomposition.

According to formulas \eqref{thipo} and \eqref{thipoodd}, the idempotent $e^{}_{\Uc}$
is the sum of primitive idempotents projecting the simple $\Sc_n$-module $U^{\la}$
to the one-dimensional subspaces spanned by the basis vectors. Therefore,
we derive from the Sergeev duality that the $\q^{}_N$-module
$M(\la)=\Ec^{}_{\Uc}(\CC^{N|N})^{\ot n}$ is isomorphic to the direct sum of
$2^{\lfloor \frac{\ell(\la)}{2}\rfloor}_{}$ copies of the simple module $L(\la)$.
Hence the first part of the following proposition is immediate from the
character formula for $L(\la)$ \cite{s:ta}; see also \cite[Theorem~3.51]{cw:dr}.
The second part was stated in \eqref{strex}.
Let $y_1,\dots,y_N$ be variables. Set $y_{-i}=y_i$ for $i=1,\dots,N$ and
introduce the diagonal matrix
\beql{y}
Y=\sum_{i=-N}^N y_i\ts E_{ii}(-1)^{\bi}.
\eeq
We will use the subscript notation for this matrix as in \eqref{tennot} and \eqref{ja}.

\bpr\label{prop:char}
The character of the $\q^{}_N$-module $M(\la)$ coincides with $Q_{\la}(y)$.
Moreover, we have the identity
\beql{topde}
\str\ts\Ec^{}_{\Uc}\tss Y_1\dots Y_n=Q_{\la}(y),
\eeq
where the supertrace is taken over all $n$ copies of the endomorphism algebra $\End\CC^{N|N}$.
\epr

\bpf
We only need to verify the second part. By the calculation in \cite[Sec.~2.2]{s:ta},
we have the identity
\ben
2^{\lfloor \frac{\ell(\la)}{2}\rfloor}_{}\str\ts\Xc^{\la}\ts Y_1\dots Y_n
=2^n\tss n!\ts Q_{\la}(y),
\een
where $\Xc^{\la}$ is the image of the character $\chi^{\la}$
under the action
of $\Sc_n$ in $(\CC^{N|N})^{\ot n}$;
see also \cite[Prop.~4.8]{n:ci}. Now apply Corollary~\ref{cor:char}
and note the relation
\beql{huhin}
\str\ts H\ts \Ec_{\Vc}\tss H^{-1}\ts Y_1\dots Y_n = \str\ts \Ec_{\Vc}\tss Y_1\dots Y_n
\eeq
which holds for the image $H$ of any basis element $h$ in \eqref{basele} under the action
of $\Sc_n$ in $(\CC^{N|N})^{\ot n}$. Indeed,
the product $Y_1\dots Y_n$ commutes with any element of $\Sym_n$
acting in $(\CC^{N|N})^{\ot n}$. Hence, writing $h$ as in \eqref{basele},
and using the cyclic property of the supertrace, we are left to verify
the relation for the basis elements of the form $h=c_1^{\ve_1}\dots c_n^{\ve_n}$.
Now observe that $J_aY_a=-Y_aJ_a$ for $a=1,\dots,n$, while $J_aY_b=Y_bJ_a$ for $a\ne b$.
Therefore, \eqref{huhin} is implied by the property that
$\str_a J_aXJ_a=\str_a X$ which holds
for an arbitrary even element $X$
of the tensor product of the endomorphism algebras.
The desired identity \eqref{topde} now follows by taking into account the number of basis elements
\eqref{basele}.
\epf

\section{Capelli identities}
\label{sec:uci}

We let $e_{kl}=E_{kl}+E_{-k,-l}$ and $f_{kl}=E_{k,-l}+E_{-k,l}$
with $k,l=1,\dots,N$,
denote the respective standard even and odd generators of the superalgebra $Q_N$.
Introduce the odd element
\beql{G}
G=\sum_{k,l=1}^N \big(e_{kl}\ot F_{l,-k}-f_{kl}\ot F_{lk}\big)
\in Q_N\ot\U(\q^{}_N).
\eeq
We will adapt a standard matrix notation to multiple tensor products
of the form
\beql{omul}
\underbrace{Q_N\ot\dots\ot Q_N}_n\ts\ot\ts \U(\q^{}_N).
\eeq
For $a=1,\dots,n$ define elements of the superalgebra \eqref{omul} by
\beql{ga}
G_a=\sum_{k,l=1}^N \Big(1^{\ot (a-1)}\ot e_{kl}\ot 1^{\ot (n-a)}\ot F_{l,-k}
-1^{\ot (a-1)}\ot f_{kl}\ot 1^{\ot (n-a)}\ot F_{lk}\Big).
\eeq
Furthermore, for $1\leqslant a<b\leqslant n$ introduce its elements
by
\begin{multline}\label{tab}
T_{ab}=\sum_{k,l=1}^N \Big(1^{\ot (a-1)}\ot f_{kl}\ot 1^{\ot (b-a-1)}\ot e_{lk}\ot 1^{\ot (n-b)}\ot 1\\
{}-1^{\ot (a-1)}\ot e_{kl}\ot 1^{\ot (b-a-1)}\ot f_{lk}\ot 1^{\ot (n-b)}\ot 1\Big).
\end{multline}
Note that
the defining relations of the superalgebra $\U(\q^{}_N)$ can now be written in a matrix form
as the relation
\beql{madef}
G_1G_2+G_2G_1=-T_{12}\tss G_1-G_1\tss T_{12}
\eeq
in the superalgebra \eqref{omul} with $n=2$.

\subsection{Odd Capelli identity}
\label{subsec:oci}

Introduce the supercommutative superalgebra $\Pc$ with generators $x_{ak}$ labelled
by $a$ running over the set $\{-M,\dots,-1,1,\dots,M\}$ for a positive integer $M$
and $k$ running over the set $\{1,\dots,N\}$. The parity of $x_{ak}$ is $\bar a\mod 2$,
where
\ben
\bar a=\begin{cases} 1\qquad\text{for}\quad a<0,\\
0\qquad\text{for}\quad a>0.
\end{cases}
\een
We let $\di_{ka}=\di/\di x_{ak}$ be the left derivation; its parity is $\bar a$.
Recall the action of $\q^{}_N$ on the space of polynomials $\Pc$,
considered in the context
of the queer Howe duality \cite[Sec.~5.2.4]{cw:dr} and \cite{s:ac}.
It is given by
\beql{howefklact}
\bal
F_{kl}&\mapsto \sum_{a=-M}^M x_{ak}\tss\di_{\tss la},\\
F_{k,-l}&\mapsto \sum_{a=-M}^M x_{ak}\tss\di_{\tss l,-a},
\eal
\eeq
with $k,l=1,\dots,N$.
We will use the composition of the action \eqref{howefklact}
with the automorphism of $\q^{}_N$ defined by
\beql{autom}
F_{kl}\mapsto -F_{lk}\Fand F_{k,-l}\mapsto i\ts F_{l,-k},\qquad k,l=1,\dots,N.
\eeq
This defines another representation
of $\q^{}_N$ in $\Pc$ given by
\beql{fklact}
\bal
F_{kl}&\mapsto -\sum_{a=-M}^M x_{al}\tss\di_{ka},\\
F_{k,-l}&\mapsto i\ts\sum_{a=-M}^M x_{al}\tss\di_{k,-a}.
\eal
\eeq
We thus get a homomorphism from
$\U(\q^{}_N)$ to the superalgebra of differential operators $\PD$
associated with $\Pc$.

To write this representation in a matrix form, introduce the auxiliary
homomorphism superspaces
$\Hom(\CC^M,\CC^{N|N})$ and $\Hom(\CC^{N|N},\CC^M)$ with the respective
sets of basis elements $\{e_{ka},f_{ka}\}$ and $\{e_{ak},f_{ak}\}$, where
$a\in\{1,\dots,M\}$ and $k\in\{1,\dots,N\}$. The basis elements $e_{ka},e_{ak}$ are even,
while $f_{ka},f_{ak}$ are odd. We have the natural product maps
\beql{hompr}
\Hom(\CC^M,\CC^{N|N})\ot\Hom(\CC^{N|N},\CC^M)\to Q_N
\eeq
defined by
\ben
e_{ka}e_{bl}=\de_{ab}e_{kl},\qquad e_{ka}f_{bl}=\de_{ab}f_{kl},\qquad f_{ka}e_{bl}=\de_{ab}f_{kl},
\qquad f_{ka}f_{bl}=\de_{ab}e_{kl},
\een
and
\ben
\Hom(\CC^{N|N},\CC^M)\ot Q_N \to \Hom(\CC^{N|N},\CC^M)
\een
with
\ben
e_{ak}e_{ml}=\de_{km}e_{al},\qquad e_{ak}f_{ml}=\de_{km}f_{al},\qquad
f_{ak}e_{ml}=\de_{km}f_{al},\qquad f_{ak}f_{ml}=\de_{km}e_{al}.
\een
Combine the variables and derivations into matrices $X$ and $D$ by setting
\ben
X=\sum_{a=1}^M\sum_{k=1}^N\big(e_{ka}\ot x_{ak}-i\ts f_{ka}\ot x_{-a,k}\big)\in
\Hom(\CC^M,\CC^{N|N})\ot \PD
\een
and
\ben
D=\sum_{a=1}^M\sum_{k=1}^N\big(i\ts e_{ak}\ot \di_{k,-a}+ f_{ak}\ot \di_{ka}\big)\in
\Hom(\CC^{N|N},\CC^M)\ot \PD,
\een
so that $X$ is even and $D$ is odd.
The representation \eqref{fklact} of $\q^{}_N$ in $\Pc$
can now be written as
\beql{gxd}
G\mapsto XD,
\eeq
where the product $XD$ is understood as an element of $Q_N\ot\PD$
defined via \eqref{hompr}.

We will
consider multiple tensor products
of the homomorphism spaces and the maps
\beql{multhom}
\big(\Hom(\CC^M,\CC^{N|N})\ot\Hom(\CC^{N|N},\CC^M)\big)^{\ot n}\to (Q_N)^{\ot n},
\eeq
where the product in each tensor factor is defined
by \eqref{hompr}. Similar to \eqref{ga}, we will write $X_r$ and $D_r$
for $r=1,\dots,n$ to denote
the matrices $X$ and $D$, associated with the $r$-th copy of the respective homomorphism space
in \eqref{multhom}.
We also need the images of the odd Jucys--Murphy elements \eqref{oddjm}
(multiplied by $\sqrt{2}$)
under the action \eqref{saact} of the Sergeev superalgebra $\Sc_n$.
They are given by
\ben
M^{(1)}=0,\qquad
M^{(b)}=T_{1b}+\dots+T_{b-1,b}\in (Q_N)^{\ot n},\qquad b=2,\dots,n,
\een
with the elements $T_{ab}$ defined in \eqref{tab}.
We will identify $M^{(b)}$ with the element $M^{(b)}\ot 1$ of the superalgebra \eqref{omul}.

We are now in a position to state the universal odd Capelli identity.

\bth\label{thm:uca}
Under the homomorphism \eqref{gxd}, we have
\ben
(G_1+M^{(1)})\dots (G_n+M^{(n)})\mapsto X_1\dots X_nD_1\dots D_n.
\een
\eth

\bpf
Arguing by induction on $n$, we come to verify the identity
\ben
X_1\dots X_{n-1}D_1\dots D_{n-1}\big(X_nD_n+T_{1n}+\dots+T_{n-1,n}\big)=
X_1\dots X_nD_1\dots D_n.
\een
We have the relations
\ben
D_rX_n=X_nD_r+\wt T_{rn}
\een
for $r=1,\dots,n-1$, where
\ben
\wt T=\sum_{a=1}^M \ts \sum_{l=1}^N\ts \big(f_{al}\ot e_{la}-e_{al}\ot f_{la}\big)
\in\Hom(\CC^{N|N},\CC^M)\ot\Hom(\CC^M,\CC^{N|N}),
\een
and the subscripts indicate the copies of $\wt T$ in the respective tensor factors in \eqref{multhom}.
Hence,
\ben
D_1\dots D_{n-1}X_n=X_nD_1\dots D_{n-1}+\sum_{r=1}^{n-1}D_1\dots \wt T_{rn}\dots D_{n-1},
\een
where $\wt T_{rn}$ takes the place of $D_r$ in the summands on the right hand side.
The desired identity now follows from the relations
\ben
\wt T_{rn}D_n+D_rT_{rn}=0,
\een
for $r=1,\dots,n-1$.
\epf

\subsection{Even Capelli identity}
\label{subsec:eci}

Regarding $Q_N$ as a subalgebra of $\End\CC^{N|N}$,
we can use \eqref{G} to consider $G$ as the element
\ben
G=\sum_{k,l=-N}^N E_{kl}\ot F_{l,-k}(-1)^{\bk+\bl}\in \End\CC^{N|N}\ot\U(\q^{}_N).
\een
In this form, it is related to the
even element
\ben
F=\sum_{k,l=-N}^N E_{kl}\ot F_{lk}(-1)^{\bl}\in \End\CC^{N|N}\ot\U(\q^{}_N)
\een
considered in \cite{n:ci} by $F=GJ$, where
$J$ is defined in \eqref{j}.

To state an even counterpart of Theorem~\ref{thm:uca}, we will
adapt the notation of Sec.~\ref{subsec:oci} to work with
the superalgebra of endomorphisms $\End\CC^{N|N}$ instead of $Q_N$
so that the tensor product superalgebra \eqref{omul} is replaced by
\ben
\underbrace{\End\CC^{N|N}\ot\dots\ot \End\CC^{N|N}}_n\ts\ot\ts \U(\q^{}_N).
\een
Its elements $F_a$ are defined by 
\beql{fa}
F_a=\sum_{k,l=-N}^N 1^{\ot (a-1)}\ot E_{kl}\ot 1^{\ot (n-a)}\ot F_{lk}(-1)^{\bl}
\eeq
for $a=1,\dots,n$.
We will need the extended homomorphism superspaces
$\Hom(\CC^{M|M},\CC^{N|N})$ and $\Hom(\CC^{N|N},\CC^{M|M})$ with the respective
standard sets of basis elements $E_{ka}$ and $E_{ak}$ of parity $\ba+\bk$, where
\ben
a\in\{-M,\dots,-1,1,\dots,M\}\Fand
k\in\{-N,\dots,-1,1,\dots,N\}.
\een
We have the natural product maps
\ben
\Hom(\CC^{M|M},\CC^{N|N})\ot\Hom(\CC^{N|N},\CC^{M|M})\to \End\CC^{N|N}
\een
defined by $E_{ka}E_{bl}=\de_{ab}E_{kl}$
and
\ben
\Hom(\CC^{N|N},\CC^{M|M})\ot \End\CC^{N|N} \to \Hom(\CC^{N|N},\CC^{M|M})
\een
with $E_{ak}E_{ml}=\de_{km}E_{al}$.
The modified matrices $X$ and $D$ are now defined by
\ben
X=\sum_{a=1}^M\sum_{k=1}^N\big((E_{ka}+E_{-k,-a})\ot x_{ak}-i\ts (E_{-k,a}+E_{k,-a})\ot x_{-a,k}\big),
\een
as an element of $\Hom(\CC^{M|M},\CC^{N|N})\ot \PD$, and
\ben
D=\sum_{a=1}^M\sum_{k=1}^N\big((E_{-a,-k}-E_{ak})\ot \di_{ka}+ i\ts (E_{-a,k}-E_{a,-k})
\ot \di_{k,-a}\big),
\een
as an element of $\Hom(\CC^{N|N},\CC^{M|M})\ot \PD$.
Note that both $X$ and $D$ are even.
The representation \eqref{fklact} can now be written
in the form
\beql{factixd}
F\mapsto X\tss D.
\eeq
Denote by $\Xc^{(b)}$ the image of the Jucys--Murphy element $x_b$ defined in \eqref{jmde},
under the action \eqref{saact} of $\Sc_n$
on the tensor product space $(\CC^{N|N})^{\ot n}$. 
Explicitly,
\ben
\Xc^{(1)}=0,\qquad \Xc^{(b)}=\sum_{a=1}^{b-1}P_{ab}(1+J_aJ_b),\qquad b=2,\dots,n.
\een
We have
the following even version of the universal Capelli identity.

\bth\label{thm:maf}
Under the action \eqref{factixd} we have
\beql{fgmxd}
\big(F_1+\Xc^{(1)}\big)\dots \big(F_n+\Xc^{(n)}\big)
\mapsto X_1\dots X_n D_1\dots D_n.
\eeq
\eth

\bpf
As in the proof of Theorem~\ref{thm:uca},
we only need to verify the identity
\ben
X_1\dots X_{n-1}D_1\dots D_{n-1}\Big(X_nD_n+\sum_{r=1}^{n-1} (P_{rn}+P_{rn}J_rJ_n)\Big)=
X_1\dots X_nD_1\dots D_n.
\een
Now we have
\ben
D_rX_n=X_nD_r-\wt T_{rn}
\een
for $r=1,\dots,n-1$, where
\ben
\wt T=\sum_{a=-M}^M \ts \sum_{l=-N}^N\ts \big(E_{al}
\ot E_{la}(-1)^{\bl}+E_{al}\ot E_{-l,-a}(-1)^{\bl}\big)
\een
is an element of
\ben
\Hom(\CC^{N|N},\CC^{M|M})\ot\Hom(\CC^{M|M},\CC^{N|N}).
\een
It remains to note that
\ben
-\wt T_{rn}D_n+D_r(P_{rn}+P_{rn}J_rJ_n)=0,
\een
for $r=1,\dots,n-1$, which is straightforward to verify.
\epf

Now fix a standard barred tableau $\Uc$ of shape $\la$
and consider the associated signed contents defined in \eqref{sgncont}.
Recall that $\Ec_{\Uc}$ denotes the image of the idempotent $e_{\Uc}$ in the tensor
product space $(\End\CC^{N|N})^{\ot n}$. The following is a queer superalgebra version
of the higher Capelli identities of \cite{o:yb}; see also \cite[Theorem~7.4.1]{m:yc}.

\bco\label{cor:oktype}
Under the action \eqref{factixd} we have
\ben
\big(F_1+\ka_1(\Uc)\big)\dots \big(F_n+\ka_{\nsf}(\Uc)\big)\ts\Ec^{}_{\Uc}
\mapsto X_1\dots X_n D_1\dots D_n\ts \Ec^{}_{\Uc}.
\een
\eco

\bpf
Multiply both expressions
appearing in \eqref{fgmxd} by $\Ec^{}_{\Uc}$ from the right.
We have the relation
\beql{fue}
\big(F_1+\Xc^{(1)}\big)\dots \big(F_n+\Xc^{(n)}\big)\Ec^{}_{\Uc}
=\big(F_1+\ka_1(\Uc)\big)\dots \big(F_n+\ka_{\nsf}(\Uc)\big)\Ec^{}_{\Uc}.
\eeq
Indeed,
let $\Vc$ by the tableau obtained from $\Uc$
by deleting the entry $\nsf$ (equal to $n$ or $\bar n$). By the definition \eqref{murphyse}
of the idempotents,
we have $e_{\Uc}=e_{\Vc}e_{\Uc}$. Due to \eqref{xietsign},
since $\Ec^{}_{\Vc}$ is permutable with $F_n$,
we can write
\ben
\big(F_n+\Xc^{(n)}\big)\Ec^{}_{\Uc}=\big(F_n+\ka_{\nsf}(\Uc)\big)\Ec^{}_{\Uc}
=\big(F_n+\ka_{\nsf}(\Uc)\big)\Ec^{}_{\Vc}\Ec^{}_{\Uc}
=\Ec^{}_{\Vc}\big(F_n+\ka_{\nsf}(\Uc)\big)\Ec^{}_{\Uc}.
\een
Then \eqref{fue} follows by
an obvious induction.
\epf

Similar to the remark we made in the Introduction concerning the equivalence
of the universal Capelli identity for $\gl_N$ to the identities in \cite{o:yb},
note that the identities of Corollary~\ref{cor:oktype} are equivalent to those
in Theorem~\ref{thm:maf}. This is clear due to the property that
the idempotents $e^{}_{\Uc}$ form a decomposition of the identity; see \cite[Prop.~2.2]{kms:jm}.

\section{Harish-Chandra images of quantum immanants}
\label{sec:hch}

Recall from \eqref{smuu} that the {\em quantum immanant} associated with a strict partition
$\la\Vdash n$ is defined by
\beql{secsmuu}
\SSb_{\la}=\str\ts \Ec^{}_{\Uc}\ts\big(F_1+\ka_1(\Uc)\big)\dots \big(F_n+\ka_{\nsf}(\Uc)\big),
\eeq
by using a fixed standard barred tableau $\Uc$ of shape $\la$
with $\ell(\la)\leqslant N$, where
the signed contents $\ka_{\asf}(\Uc)$ are
defined in \eqref{sgncont}. We will see below in Theorem~\ref{thm:hch}, that the quantum immanant
depends only on $\la$ and does not depend on $\Uc$.

We keep the notation $\Zr(\q^{}_N)$ for the center of the universal enveloping algebra $\U(\q^{}_N)$.

\bpr\label{prop:cent}
The quantum immanant $\SSb_{\la}$ belongs to $\Zr(\q^{}_N)$.
\epr

\bpf
Consider the tensor product superalgebra
\beql{addit}
(\End\CC^{N|N})^{\ot (n+1)}\ot \U(\q^{}_N)
\eeq
with an additional copy of $\End\CC^{N|N}$ labelled by $n+1$. It will be sufficient to
verify that $\SSb_{\la}$, regarded as an element of \eqref{addit}
with the identity component in the additional copy of $\End\CC^{N|N}$,
commutes with $F_{n+1}$ in this superalgebra.
The defining relations of $\U(\q^{}_N)$
are well-known to have
the matrix form
\ben
[F_1,F_2]=(P_{12}+P_{12}J_1J_2)\tss F_1-F_1\tss (P_{12}+P_{12}J_1J_2)
\een
which can be checked directly or derived from
\eqref{madef} by multiplying both sides by $J_1J_2$ from the right. Hence,
\beql{commuq}
[F_{n+1},F_r]=F_r(P_{r,n+1}+P_{r,n+1}J_rJ_{n+1})-(P_{r,n+1}+P_{r,n+1}J_rJ_{n+1})F_r
\eeq
for $r=1,\dots,n$. Therefore,
\begin{multline}
\big[F_{n+1},\ts \Ec^{}_{\Uc}\ts\big(F_1+\ka_1(\Uc)\big)\dots \big(F_n+\ka_{\nsf}(\Uc)\big)\big]\\
{}=\sum_{r=1}^n \Ec^{}_{\Uc}\ts\big(F_1+\ka_1(\Uc)\big)\dots
[F_{n+1},F_r]\dots
\big(F_n+\ka_{\nsf}(\Uc)\big).
\non
\end{multline}
Applying \eqref{commuq}, we can write this expression as
\ben
\Ec^{}_{\Uc}\ts\big(F_1+\ka_1(\Uc)\big)\dots \big(F_n+\ka_{\nsf}(\Uc)\big) \Xc^{(n+1)}
-\Ec^{}_{\Uc}\ts\Xc^{(n+1)}\ts\big(F_1+\ka_1(\Uc)\big)\dots \big(F_n+\ka_{\nsf}(\Uc)\big),
\een
where
\ben
\Xc^{(n+1)}=\sum_{r=1}^n (P_{r,n+1}+P_{r,n+1}J_rJ_{n+1}).
\een
Note that $\Xc^{(n+1)}$ coincides with the image of the Jucys--Murphy element $x_{n+1}\in\Sc_{n+1}$
in the tensor product superalgebra and so commutes with $\Ec^{}_{\Uc}$. Hence,
swapping $\Ec^{}_{\Uc}$ and $\Xc^{(n+1)}$ in the second term in
the above expression and taking the supertrace, we find that the two terms cancel
by the cyclic property of supertrace.
\epf

The center $\Zr(\q^{}_N)$ is isomorphic to the algebra
of supersymmetric polynomials $\Ga_N$ which we recalled in Sec.~\ref{sec:sqs},
via the {\em Harish-Chandra isomorphism}
\beql{hchiso}
\chi:\Zr(\q^{}_N)\to \Ga_N;
\eeq
see \cite[Sec.~2.3]{cw:dr} and \cite{s:ip}. In the definition of $\chi$ we use the standard
(upper-triangular) Borel subalgebra spanned by the basis elements $F_{ij}$ with $|i|\leqslant |j|$
as in \cite[Sec.~2.3.2]{cw:dr}.
The following is our final principal result where we use the factorial Schur $Q$-polynomials;
see \eqref{qstar}.

\bth\label{thm:hch}
Under the Harish-Chandra isomorphism we have
\ben
\chi:\SSb_{\la}\mapsto Q^+_{\la}(y).
\een
In particular, $\SSb_{\la}$ depends only on the shape of $\Uc$.
Moreover, the quantum immanants $\SSb_{\la}$ with $\ell(\la)\leqslant N$
form a basis of $\Zr(\q^{}_N)$.
\eth

\bpf
The image $\chi(\SSb_{\la})$ is a supersymmetric polynomial in $y$ of degree
not exceeding $n$. We will use the
characterization property of the factorial Schur $Q$-polynomials
recalled in Sec.~\ref{sec:sqs} to show that $\chi(\SSb_{\la})=Q^+_{\la}(y)$.

The quantum immanant $\SSb_{\la}$ acts by scalar multiplication
in any highest weight $\q^{}_N$-module $L(y)$
with the highest weight $y=(y_1,\dots,y_N)$
with respect to the standard Borel subalgebra; here
we regard the components $y_i$ as complex numbers.
As the components $y_i$ vary, the scalar, regarded as a function of $y$,
is a supersymmetric polynomial, and it coincides with $\chi(\SSb_{\la})$.

The highest vector of $L(y)$ is an eigenvector for
the operators
$F_{11},\dots,F_{NN}$ with the respective eigenvalues
$y_1,\dots,y_N$. Therefore, the top degree component
of the polynomial $\chi(\SSb_{\la})$ can be written as
the supertrace
$
\str\ts\Ec^{}_{\Uc}\ts Y_1\dots Y_n
$
of the diagonal matrix defined in \eqref{y}. Hence,
by \eqref{topde},
the top degree component
of $\chi(\SSb_{\la})$ coincides with $Q_{\la}(y)$.

To apply the
characterization property of the factorial Schur $Q$-polynomials,
we only need to show that $\chi(\SSb_{\la})$ vanishes
when $y$ is evaluated as $y=-\mu$ for all shifted Young diagrams $\mu$
with $|\mu|<n$. It will be sufficient to show that
the central element $\chi(\SSb_{\la})$ acts as zero in
the corresponding highest weight modules $L(-\wt\mu)$, where
$\wt\mu=(\mu_N,\dots,\mu_1)$.
To this end, we will use
the realization of $L(\mu)$ in the space of polynomials $\Pc$
with the action of $\q^{}_N$ given by \eqref{howefklact}.
We will assume that
the parameter $M$ is large enough (at least, $M\geqslant N$). The highest weight module
$L(\mu)$ is then realized as a submodule of $\Pc$, and consists of homogeneous
polynomials in the $x_{ak}$
of degree $|\mu|$. The
action
\eqref{factixd} is obtained from \eqref{howefklact} by twisting by
the automorphism \eqref{autom}. Under the twisted action, $L(\mu)$ becomes
the highest weight module with the highest weight $-\mu$ with respect to the
Borel subalgebra which is opposite to the standard upper-triangular subalgebra.
This module is isomorphic to the dual module $L(\mu)^*$, and hence
is isomorphic to the highest weight module $L(-\wt\mu)$ with respect to the standard Borel
subalgebra; see e.g. \cite[Exercise~2.18]{cw:dr}.

On the other hand, by Corollary~\ref{cor:oktype}, the image of
$\SSb_{\la}$ in the superalgebra $\PD$ is a differential operator
which annihilates all polynomials in the $x_{ak}$ of degree not exceeding $n-1$.
Thus, $\SSb_{\la}$ acts as zero in the modules $L(-\wt\mu)$ with $|\mu|<|\la|=n$,
as required.

The last part of the theorem holds due to the isomorphism \eqref{hchiso},
since the polynomials $Q^+_{\la}(y)$ form a basis of $\Ga_N$.
\epf

We keep the notation $\Xc^{\la}$ for the image of the character $\chi^{\la}$
under the action
of $\Sc_n$ in $(\CC^{N|N})^{\ot n}$.
The following is an analogue of the Capelli identities
of \cite[Theorem~1.3]{ass:sq} and \cite[Theorem~4.7]{n:ci}.

\bco\label{cor:inde}
The image of $\SSb_{\la}$ under the action \eqref{factixd} can be written in the form
\ben
\SSb_{\la}\mapsto \frac{2^{\lfloor \frac{\ell(\la)}{2}\rfloor}_{}}{2^n\ts n!}\ts
\str\ts \Xc^{\la} X_1\dots X_n D_1\dots D_n.
\een
\eco

\bpf
This follows from Corollary~\ref{cor:oktype} by the same argument as in the proof of
Proposition~\ref{prop:char}. In the notation of that proof,
it suffices to verify the identity
\ben
\str\ts H\ts \Ec_{\Vc}\tss H^{-1}\ts X_1\dots X_n D_1\dots D_n
= \str\ts \Ec_{\Vc}\tss X_1\dots X_n D_1\dots D_n
\een
which is a counterpart of \eqref{huhin}.
The product $X_1\dots X_n D_1\dots D_n$ commutes with any element of $\Sym_n$
acting in $(\CC^{N|N})^{\ot n}$, while
\ben
J_aX_1\dots X_n D_1\dots D_n=-X_1\dots X_n D_1\dots D_n J_a
\een
for $a=1,\dots,n$. These properties are seen by introducing
the elements
\ben
\wt P=\sum_{a,b=-M}^M E_{ab}\ot E_{ba}(-1)^{\bb}\in\End\CC^{M|M}\ot \End\CC^{M|M}
\een
and
\ben
\wt J=\sum_{a=-M}^M E_{a,-a}(-1)^{\ba}\in \End\CC^{M|M}
\een
and checking the relations
\ben
P_{12}\ts X_1X_2=X_1X_2\ts\wt P_{12}\Fand  \wt P_{12}D_1D_2=D_1D_2 P_{12}
\een
along with
\ben
JX=X\wt J \Fand \wt J D=-DJ.
\een
The final steps are the same as in the proof of Proposition~\ref{prop:char}.
\epf

\bre\label{rem:concz}
As a concluding remark, we point out a relation of our formulas with
those of \cite{ass:sq} and \cite{n:ci}.
The element $\Psi_{\la}\in\Sc_n$ is defined in \cite[Sec~2]{n:ci}
via a fusion procedure. According to
its version obtained in \cite[Theorem~4.1]{kms:jm}, we find that
\ben
\Psi_{\la}=\frac{n!}{g_{\la}}\ts e^{}_{\Uc^c}
\een
for the standard tableau $\Uc^c$ obtained by filling the boxes of $\la$
with the numbers $1,\dots,n$
by successive columns from top to bottom. Therefore, by extending the argument
of \cite[Prop.~3.4]{o:qi} to the Sergeev superalgebra, we can derive a relation
for the Capelli element $C_{\la}$ of \cite{n:ci}
in terms of an arbitrary standard barred tableau $\Uc$ of shape $\la$:
\ben
\frac{g_{\la}}{n!}\ts C_{\la}=\str\ts \Ec^{}_{\Uc}\ts\big(F_1-\ka_1(\Uc)\big)\dots \big(F_n-\ka_{\nsf}(\Uc)\big);
\een
cf.
\eqref{secsmuu}.
By the same argument as in the proof of Theorem~\ref{thm:hch}, relying on the Capelli
identity of \cite[Theorem~4.7]{n:ci}, we get a formula for the Harish-Chandra
image of $C_{\la}$:
\beql{caim}
\frac{g_{\la}}{n!}\ts \chi(C_{\la})=Q^-_{\la}(y),
\eeq
as defined in \eqref{qstar};
cf. \cite[Corollary~1.5]{ass:sq}.
We can derive from the proof of Theorem~\ref{thm:hch}
that $\frac{(-1)^n\tss g_{\la}}{n!}\ts C_{\la}$
coincides with the image of the quantum immanant $\SSb_{\la}$
under the anti-automorphism of the superalgebra $\U(\q^{}_N)$ taking $F_{kl}$ to $-F_{kl}$.

According to \cite[Theorem~1.3]{ass:sq}, the Harish-Chandra
image of the central element $z_{\la}$ produced therein, is given by\footref{sahi}
\ben
\chi(z_{\la})=\frac{(-1)^{n}\ts g_{\la}}{2^{\tss\ell(\la)}_{}\ts n!}\ts Q^+_{\la}(y),
\een
and so
\beql{zlasla}
z_{\la}= \frac{(-1)^{n}\ts g_{\la}}{2^{\tss\ell(\la)}_{}\ts n!}\ts \SSb_{\la}
\eeq
by Theorem~\ref{thm:hch}.
\ere

\section*{Declarations}

\subsection*{Competing interests}
The authors have no competing interests to declare that are relevant to the content of this article.

\subsection*{Acknowledgements}
We are grateful to Siddhartha Sahi and Hadi Salmasian for a discussion
of the results of \cite{ass:sq}.
The work on
the project was completed during the second author's visit
to the Shenzhen International Center for Mathematics. He
is grateful to the Center for warm hospitality. His
work was also supported by the Australian Research Council, grant DP240101572.

\subsection*{Availability of data and materials}
No data was used for the research described in the article.

\newpage


\small

\noindent
I.K.:\newline
Shenzhen International Center for Mathematics\\
Southern University of Science and Technology, China\\
{\tt kashuba@sustech.edu.cn}

\vspace{5 mm}

\noindent
A.M.:\newline
School of Mathematics and Statistics\newline
University of Sydney,
NSW 2006, Australia\newline
{\tt alexander.molev@sydney.edu.au}

\end{document}